\documentstyle[11pt]{article}

\include{amssym}
\newcommand{\R}{\makebox{$\Bbb R$}}
\newcommand{\dpl}{\displaystyle}
\newcommand{\qed}{\raisebox{-.15ex}{\rule{1.5mm}{1.5ex}}}
\begin{document}
\noindent \hphantom{a} \hfill {\Large\sf{\textbf{Some incomplete and
boundedly complete}}}
\hfill \hphantom{a} \\
\hphantom{a} \hfill {\Large\sf{\textbf{families of discrete
distributions}}}\footnote{\small\sf{\textbf{TECHNICAL REPORT No.
BIRU/2009/2; Date: September 17, 2009}}} \hfill \hphantom{a}

\vspace{0.17in}\noindent \centerline{\Large\sf Sumitra Purkayastha}

\vspace{0.05in}\noindent \centerline{\large\em Bayesian and
Interdisciplinary Research Unit} \centerline{\large\em Indian
Statistical Institute} \centerline{\large\em 203 B.T. Road. Kolkata
700 108. India} \centerline{\large\sf E-mail: sumitra@isical.ac.in}

\vspace{0.17in} \noindent \centerline{\large\it {\textbf{Abstract}}}

%\vspace{0.1in}
\noindent \begin{quote} {\em We present a general result giving us
families of incomplete and boundedly complete families of discrete
distributions. For such families, the classes of unbiased estimators
of zero with finite variance and of parametric functions which will
have uniformly minimum variance unbiased estimators with finite
variance are explicitly characterized. The general result allows us
to construct a large number of families of incomplete and boundedly
complete families of discrete distributions. Several new examples of
such families are described.}
\end{quote}

\vspace{0.05in} \noindent {\sf {\textbf{Keywords}}} Boundedly
complete family of distribution, complete family of distribution,
discrete distribution, power series distribution
\setcounter{page}{1}

\vspace{0.17in} \noindent {\large\sf {\textbf{1. Introduction}}}

\vspace{0.05in}\noindent In Lehmann and Scheffe (1950, p. 312)
(available also in Lehmann and Casella (1998, chapter 2, section 1,
pp. 84--85)), an example of a boundedly complete but not complete
family of distributions was given. This family has the following
property: not every unbiasedly estimable parametric function will
have a uniformly minimum variance unbiased estimator (UMVUE). Also,
for this family, the classes of unbiased estimators of zero with
finite variance and of parametric functions which will have UMVUE's
with finite variance are explicitly characterized. On the other
hand, not many examples of families of distributions, especially
{\bf discrete}, boundedly complete but not complete, are available
in the literature. This had been noticed earlier and necessary
investigations were accordingly undertaken by several authors. Such
an investigation is meaningful not only for its own sake but also
for introducing completeness and explaining its role in
understanding unbiased estimators. Earlier investigations in this
direction include Hoeffding (1977), Bar-Lev and Plachky (1989), and
Mattner (1993).

Hoeffding's (1977) examples of incomplete and boundedly complete
families of distributions were in non-parametric set-up. Bar-Lev and
Plachky (1989) obtained similar examples of discrete distributions.
Mattner's (1993) examples were obtained for location families of
probability measures on Euclidean space.

Taking an approach, different from that in Bar-Lev and Plachky
(1989), we present in this paper a general result giving us families
of discrete distributions which are boundedly complete but not
complete. Characterizations of unbiased estimators of zero with
finite variance and the parametric functions which will have UMVUE's
with finite variance are also parts of the result. This result is
stated and proved in section 2. Based on this result, we prove in
section 3 a general result that enables us to construct from {\bf
any} given power series distribution incomplete and boundedly
complete families of discrete distributions. This result involves
choice of a suitable infinite sequence of real numbers $\{b_k: k
\geq 0\}$. Accordingly, in section 3 we also identify {\bf specific}
choices of $\{b_k: k \geq 0\}$ which will work for {\bf any} given
power series distribution. Some new examples of incomplete and
boundedly complete families related to {\bf specific} power series
distributions are given in section 4.

\vspace{0.17in} \noindent {\large\sf {\textbf{2. A general result}}}

\vspace{0.05in} \noindent We prove a general result in this section
that allows us to construct from a given power series distribution,
generated by a power series $f$, families of discrete distributions
which are boundedly complete but not complete. We characterize also
the corresponding classes of unbiased estimators of zero with finite
variance and of parametric functions which will have UMVUE's with
finite variance.

\vspace{0.05in} \noindent {\sf {\textbf{Theorem 1}}} {\em Consider a
power series $f(\theta)$ given by} $f(\theta) =
\sum_{k=0}^{\infty}a_k\theta^k,$ {\em where $a_k > 0\;\forall\;k
\geq 0$, and the radius of convergence $R_1$ is positive. Assume the
following.}
\begin{enumerate}
\item[(A)] {\em  There exists an unbounded sequence of numbers $\{b_k: k \geq 0\}$
such that $b_0 = 1$, $b_k \geq 1\;\forall\;x \geq 1$, the set $J :=
\{k \geq 0: b_k =1\}$ is finite, and the power series $h(\theta) :=
\sum_{k=0}^{\infty} b_k^2 a_k \theta^k$ has a positive number $R_2$
as its radius of convergence.}
\end{enumerate} {\em Let $R := \min(R_1,R_2).$ ($R$ may be $+\infty$.)
Define $g(\theta) := \sum_{k=0}^{\infty} b_k a_k \theta^k$ for
$\theta \in (0,R)$. Consider the family ${\cal P}_{f,g}$ of discrete
distributions given by} ${\cal P}_{f,g} = \{P_{\theta}: \theta \in
(0,R)\}$, {\em where the probability mass function (pmf)
$p(k;\theta)$ corresponding to $P_{\theta}$ is given by}
$$p(k;\theta) = \left\{\begin{array}{lcl} 1- f(\theta)/g(\theta) & {\mbox{\it if}} & k = -1, \\
a_k \theta^k/g(\theta), & {\mbox{\it if}} & k = 0,1,2,\ldots
.\end{array}\right.$$ {\em Suppose $X$ has pmf $p(k;\theta),\;\theta
\in (0,R)$. Then the following hold.}
\begin{enumerate}
\item[(a)] {\em $\delta(X)$ is an unbiased
estimator of zero with finite variance if and only if $$\delta(k) =
-\delta(-1)\left(b_k-1\right)\;{\mbox{\it for}}\;k = 0, 1,2,\ldots.
\eqno{(1)}$$}
\item[(b)] {\em ${\cal P}_{f,g}$ is not complete but boundedly complete.}

\item[(c)] {\em A parametric function $\psi(\theta)$ has a
UMVUE with finite variance if and only if $\psi(\theta) = a +
\sum_{k\in J} c_k \theta^k/g(\theta)$ for some $a, c_k\;(k \in J)\;
\in {\R}$.}
\end{enumerate}

\noindent {\sf {\textbf{Remark 1}}} From $(\sum_{k=0}^{\infty} b_k
a_k \theta^k)^2 \leq (\sum_{k=0}^{\infty} b_k^2 a_k \theta^k)\cdot
(\sum_{k=0}^{\infty} a_k \theta^k)$ for any $\theta > 0$, it can be
seen that convergence of $\sum_{k=0}^{\infty} b_k^2 a_k \theta^k$
implies that of $\sum_{k=0}^{\infty} b_k a_k \theta^k$, if
$\sum_{k=0}^{\infty} a_k \theta^k$ is convergent. In other words,
$g$ is well-defined.

\vspace{0.05in} \noindent {\sf {\textbf{Remark 2}}} Condition (A)
implies that $b_k >1$ for all $k$ sufficiently large. Hence, we get
the following: for every $\theta \in (0,R)$, $0 < f(\theta) <
g(\theta)$. Consequently, $0 < p(-1;\theta) < 1$ and $p(k;\theta)$
is well-defined.

\vspace{0.05in} \noindent {\sf {\textbf{Remark 3}}} In the sequel,
we shall call ${\cal P}_{f,g}$, {\em the family of discrete
distributions induced by $f$ and $g$.} Also, we denote by ${\cal
U}_0$, the class of unbiased estimators $\delta(X)$ of zero with
finite variance. Also, we denote by $\Psi_0$, the class of functions
in part (c) of the theorem.

\vspace{0.05in} \noindent {\sf {\textbf{Proof of Theorem 1}}}

(a) To begin with, note that if (1) is satisfied, it follows from
condition (A) that $\sum_{k=0}^{\infty}b_k^2 a_k\theta^k p(k;\theta)
< \infty$ $\forall$ $\theta \in (0,R)$, implying ${\mbox{\rm
E}}_{\theta}(\delta^2(X)) < \infty$ $\forall$ $\theta \in (0,R)$.
Notice now that ${\mbox{\rm E}}_{\theta}(\delta(X)) =
\delta(-1)\left\{1- f(\theta)/g(\theta)\right\} +
[g(\theta)]^{-1}\sum_{k=0}^{\infty}\delta(k)a_k\theta^k. $
Therefore, ${\mbox{\rm E}}_{\theta}(\delta(X))$ $=0$ if and only if
$\sum_{k=0}^{\infty}\delta(k)a_k\theta^k =
-\delta(-1)\{g(\theta)-f(\theta)\}.$ Observe that
$g(\theta)-f(\theta) = \sum_{k=0}^{\infty}a_k(b_k-1)\theta^k$.
Hence, ${\mbox{\rm E}}_{\theta}(\delta(X)) =0\;\forall\;\theta \in
(0,R)$ if and only if (1) is satisfied. \hfill \qed

(b) Choose in (1), $\delta(-1) = -1$. Note that with such a choice
of $\delta(\cdot)$, $P_{\theta}(\delta(X) \neq 0) \geq P_{\theta} (X
= -1) = 1- f(\theta)/g(\theta) > 0$, for any $\theta \in (0,R)$. On
the other hand, since $\delta(\cdot)$ satisfies (1), ${\mbox{\rm
E}}_{\theta}(\delta(X)) = 0\;\forall\;\theta \in (0,R)$. Hence,
${\cal P}_{f,g}$ cannot be complete.

To see that ${\cal P}_{f,g}$ is boundedly complete, suppose
${\mbox{\rm E}}_{\theta}(\delta(X)) = 0\;\forall\;\theta \in (0,R)$
for some bounded function $\delta(\cdot)$. Note, in view of fact
(a), that $\delta(\cdot)$ satisfies (1). Since the sequence
$\left\{b_k: k \geq 0\right\}$ is unbounded and also since
$\delta(\cdot)$ is assumed to be a bounded function, we must have
$\delta(-1) = 0$. This implies, in view of (1), $\delta(k) =
0\;{\mbox{\rm for}}\;k = -1,0, 1,2,\ldots.$ Hence, ${\cal P}_{f,g}$
is boundedly complete. \hfill \qed

(c) First suppose $\psi(\theta) = a + \sum_{k\in J} c_k
\theta^k/g(\theta)$ for some $a, c_k\;(k \in J)\; \in {\R}.$ Define
$T(X) = a+ \sum_{k\in J} (c_k/a_k) \cdot 1_{\{X = k\}}.$ Note that
${\mbox{\rm E}}_{\theta}(T(X)) =\psi(\theta)$ and ${\mbox{\rm
E}}_{\theta}(T^{2}(X)) < \infty$ $\forall\;\theta \in (0,R)$. Let
${\cal U}_0 := \{\delta(\cdot): {\mbox{\rm E}}_{\theta}(\delta(X))
=0, {\mbox{\rm E}}_{\theta}(\delta^{2}(X)) < \infty\;\forall\;\theta
\in (0,R)\}$. In other words, ${\cal U}_0$ denotes the class of
unbiased estimators $\delta(X)$ of zero with finite variance.
Observe now, in view of part (a), that for any $\delta(\cdot) \in
{\cal U}_0$, $\delta(k) = 0$ for any $k \in J$, as $b_k = 1$ for all
$k \in J$. So, for any $k \in J$, the random variable $1_{\{X = k\}}
\cdot \delta(X)$ is identically zero. Hence, it has expectation
zero. Consequently, for $\theta \in (0,R)$, ${\mbox{\rm
E}}_{\theta}(\delta(X) \cdot T(X)) = a\,{\mbox{\rm
E}}_{\theta}(\delta(X)) + \sum_{k\in J} (c_k/a_k)\,{\mbox{\rm
E}}_{\theta}(1_{\{X = k\}} \cdot \delta(X))$ $= a \cdot 0 +
\sum_{k\in J} (c_k/a_k)\cdot 0 = 0.$ Therefore, in view of theorem
1.7 (chapter 2, section 1, Lehmann and Casella (1998)), $T(X)$ is
the UMVUE of $\psi(\theta)$.

Conversely, suppose $\psi(\theta)$ has a UMVUE with finite variance.
Let $T(X)$ be this UMVUE. Then, in view of theorem 1.7 (chapter 2,
section 1, Lehmann and Casella (1998)), for any unbiased estimator
$\delta(\cdot) \in {\cal U}_0$, ${\mbox{\rm E}}_{\theta}(T(X) \cdot
\delta(X)) = 0$ $\forall\;\theta \in (0,R).$ In other words, for any
such $\delta(\cdot)$, $T(X) \cdot \delta(X) \in {\cal U}_0$. Hence,
in view of fact (a) above, we get the following:
$$T(k) \cdot \delta(k) = -T(-1) \cdot
\delta(-1)\left(b_k-1\right)\;{\mbox{\rm for}}\;k = 0,1,2,\ldots
.\eqno{(2)}$$ Observe now that $\delta(\cdot)$ satisfies (1), as
$\delta(\cdot) \in {\cal U}_0$. Therefore, the right-hand side of
(2) equals $T(-1) \cdot \delta(k)$. Consider now $\delta(\cdot)$
such that $\delta(-1) \neq 0$, and hence in view of (2), (1), and
condition (A) in the statement of the theorem, $\delta(k) \neq 0$
for $k \not \in J$. Hence, dividing both sides of (2) by $\delta(k)
\neq 0$, we get $T(k) = T(-1)$ $\forall\;k \not \in J$. In other
words, $T(k) = T(-1) + \sum_{y\in J} \{T(y)-T(-1)\}\cdot 1_{\{k =
y\}}.$ Therefore, $\psi(\theta) = a + \sum_{k\in J} c_k
\theta^k/g(\theta)$ for some $a, c_k\;(k \in J)\; \in {\R}$, where
$a := T(-1)$, $c_k := a_k(T(k)-T(-1))$. \hfill \qed

\vspace{0.05in}\noindent {\sf {\textbf{Remark 4}}} It can be seen
that $-1 \in$ support of $X$ can be replaced by any real number
$\not \in \{0, 1, 2, \ldots\}$.

\vspace{0.17in} \noindent {\large\sf {\textbf{3. Examples for
arbitrary power series distributions}}}

\vspace{0.05in} \noindent Theorem 1 in the preceding section
involves choice of the sequence $\{b_k: x = 0,1, \ldots\}$. We prove
in this section a general result (theorem 2) that gives us
conditions on $\{b_k: x = 0,1, \ldots\}$ so that sequences
satisfying these conditions will work for {\bf any} $f$. This result
is used to identify {\bf specific} choices of $\{b_k: x = 0,1,
\ldots\}$ which will work for {\bf any} $f$. Such choices are
prescribed in corollaries 1--4.

\vspace{0.05in} \noindent {\sf {\textbf{Theorem 2}}} {\em Let
$f(\cdot)$, $a_k$'s, $R_1$ be as in Theorem 1. Suppose there exists
an unbounded sequence of numbers $\{b_k: k \geq 0\}$ such that $b_0
= 1$, $b_k \geq 1\;\forall\;k \geq 1$, $\{k \geq 0: b_k =1\}$ is
finite, and $\lim_{k \rightarrow\infty} b_{k}^{1/k} = L$, for some
finite $L > 0$. Then the radius of convergence, denoted by $R_2$, of
the power series $h(\theta) := \sum_{k=0}^{\infty} b_k^2 a_k
\theta^k$ is given by $R_2 = R_1/L^2$. Consequently, $\{b_k: k \geq
0\}$ satisfies condition (A) of theorem 1.}

\vspace{0.05in} \noindent {\sf {\textbf{Proof}}} Notice that we need
only to show that $R_2 = R_1/L^2$. This is implied by (a) $R_{1} =
1/\limsup_{x\rightarrow\infty} |a_k|^{1/x}$, (b) $R_{2} =
1/\limsup_{x\rightarrow\infty} |b_k^2 a_k|^{1/x}$, and (c) $\lim_{k
\rightarrow\infty} b_{k}^{1/k} = L >0$. \hfill \qed

\vspace{0.05in}\noindent {\sf {\textbf{Remark 5}}} We may replace
the condition \lq\lq $\lim_{k \rightarrow\infty} b_{k}^{1/k} = L$,
for some finite $L > 0$\rq\rq \ by the following stronger one:
\lq\lq $\lim_{k \rightarrow\infty} b_{k+1}/b_{k} = L$, for some
finite $L > 0$\rq\rq.

\vspace{0.05in} In corollaries 1-4 below, the sequence $\{b_k: k
\geq 0\}$ can easily be seen to satisfy $b_0 = 1$, $b_k \geq
1\;\forall\;k \geq 0$. Also, it is unbounded and $\{k \geq 0: b_k
=1\}$ is finite, as $\lim_{k \rightarrow\infty} b_{k} = \infty$. We
prove this last fact only for corollary 3, the proofs corresponding
to corollaries 1, 2, and 4 being trivial. Also, we prove in each
corollary that $\lim_{k \rightarrow\infty} b_{k}^{1/k} = L$, for
some finite $L > 0$. In view of theorem 2, these will show that for
each of the corollaries 1--4, $\{b_k: k \geq 0\}$ satisfies
condition (A) of theorem 1 for any $f$ there. Let us state here that
{\em in each of the corollaries and in theorem 3, $w_j$'s are
positive numbers with $\sum_{j=1}^{n}w_j \geq 1$.}

\vspace{0.05in}\noindent {\sf {\textbf{Corollary 1}}} {\em Let
$b_{0} := 1, b_k
:= \sum_{j=1}^{n}w_j d_{j,k}^{k}$, $k \geq 1$, where %$w_j$'s are
%positive numbers with $\sum_{j=1}^{n}w_j \geq 1$,
$d_{j,k} \geq 1\;\forall\;j=1, \ldots ,n,\;k \geq 1$, and for every
$j =1, \ldots ,n$, $\lim_{k\rightarrow\infty}d_{j,k} = d_j$ exists
with $d_j > 1$. Then $\lim_{k \rightarrow\infty} b_{k}^{1/k} =
\max_{1\leq j \leq n} d_j$.}

\vspace{0.05in} \noindent {\sf {\textbf{Proof}}} Observe that
$w_{i_0}d_{i_0,k}^{k} \leq b_k \leq w d_{k}^{k}$, where $w :=
\max_{1\leq j \leq n} w_j$, $i_0$ is such that $d_{i_0} =
\max_{1\leq j \leq n} d_j$, and $d_k := \max_{1\leq j \leq n}
d_{j,k}$. Hence, $w_{i_0}^{1/k}d_{i_0,k}\leq b_k^{1/k} \leq w^{1/k}
d_{k}$. Notice now that $\lim_{k\rightarrow\infty}d_k = \max_{1\leq
j \leq n} d_j = d_{i_0}$. Hence, $\lim_{k \rightarrow\infty}
b_{k}^{1/k} = d_{i_0}$. \hfill \qed

\vspace{0.05in}\noindent {\sf {\textbf{Corollary 2}}} {\em Let
$b_{0} := 1, b_k := \sum_{j=1}^{n}w_j k^{p_j}$, $k \geq 1$, where
$p_j > 0\;\forall\;j=1, \ldots ,n.$ Then $\lim_{k \rightarrow\infty}
b_{k}^{1/k} = 1$.}

\vspace{0.05in} \noindent {\sf {\textbf{Proof}}} Observe that $1
\leq b_k \leq w k^{p^{\ast}}$, where $w := \max_{1\leq j \leq n}
w_j$ and $p^{\ast} := \max_{1\leq j \leq n} p_j$. Hence, $\lim_{k
\rightarrow\infty} b_{k}^{1/k} = 1$. \hfill \qed

\vspace{0.05in}We need two propositions, stated in the appendix, in
the next corollary.

\vspace{0.05in}\noindent {\sf {\textbf{Corollary 3}}} {\em Let
$b_{0} := 1, b_k := \sum_{j=1}^{n}w_j t(\beta_j,k)/t(\alpha_j,k)$,
$k \geq 1,$ where $0 < \alpha_j < \beta_j\;\forall\;j=1, \ldots ,n$
and $t(\gamma,k)$ is defined in proposition 1. Then $\lim_{k
\rightarrow\infty} b_{k}^{1/k} = 1$.}

\vspace{0.05in} \noindent {\sf {\textbf{Proof}}} For $\alpha, \beta$
satisfying $0 < \alpha < \beta$, let $u_k :=
t(\beta,k)/t(\alpha,k),\;k \geq 0.$ Observe that $u_0 = 1$,
$u_{k+1}/u_{k} = (\beta+k)/(\alpha+k) >1$ for $k \geq 0$, and in
view of proposition 2, $\lim_{k\rightarrow\infty} u_k = \infty$.
Hence, the sequence $\{b_k: k \geq 0\}$ is unbounded, $b_0 = 1$,
$b_k \geq 1\;\forall\;k \geq 1$, and $\{k \geq 0: b_k =1\}$ is
finite.

Observe next that $u_{k+1}/u_{k} = (\beta+k)/(\alpha+k)$ for $k \geq
0$ implies  $\lim_{k\rightarrow\infty} u_{k+1}/u_{k} = 1$ which, in
turn, implies $\lim_{k\rightarrow\infty} u_k^{1/k} = 1$. Notice now
that $1 \leq b_k \leq w t(\beta^{\ast},k)/t(\alpha_{\ast},k)$, where
$w := \max_{1\leq j \leq n} w_j$, $\alpha_{\ast} := \min_{1\leq j
\leq n} \alpha_j$, and $\beta^{\ast} := \max_{1\leq j \leq n}
\beta_j$. This fact and the preceding arguments imply that $\lim_{k
\rightarrow\infty} b_{k}^{1/k} = 1$. \hfill \qed

\vspace{0.05in}\noindent {\sf {\textbf{Corollary 4}}} {\em Let
$b_{0} := 1, b_k := \sum_{j=1}^{n}w_j \log p_{j}(k)$, $k \geq 1,$
where for every $j=1, \ldots ,n$, $p_{j}(\cdot)$ is a non-constant
polynomial with all coefficients non-negative and $p_{j}(1) \geq e$.
Then $\lim_{k \rightarrow\infty} b_{k}^{1/k} = 1$.}

\vspace{0.05in} \noindent {\sf {\textbf{Proof}}} Observe that for
any non-constant polynomial $p(\cdot)$, with all coefficients
non-negative, $\lim_{k\rightarrow \infty} p(k)$ $=\infty$,
$\lim_{k\rightarrow \infty} [p(k+1)/p(k)] =1$. Hence,
$\lim_{k\rightarrow \infty} [\log p(k+1)/ \log p(k)] =1$ which, in
turn, implies $\lim_{k\rightarrow \infty} [\log p(k)]^{1/k} = 1$.
Notice now that $1 \leq b_k \leq w \log p_{0}(k)$, where $w :=
\max_{1\leq j \leq n} w_j$ and $p_{0}(k) := \sum_{j=1}^{n}w_j
p_{j}(k)$. Also, $p_{0}(k)$ is a non-constant polynomial, with all
coefficients non-negative. This fact and the preceding arguments
imply that $\lim_{k \rightarrow\infty} b_{k}^{1/k} = 1$. \hfill \qed

\vspace{0.05in}A close scrutiny of corollaries 2--4 reveals that it
is possible to combine them by considering appropriate linear
combinations of sequences $\{b_k: k \geq 0\}$, chosen from the class
of all possible sequences suggested in them. Keeping this in mind,
we let ${\cal B}_{1}$ denote the class of all possible sequences
$\{b_k: k \geq 0\}$ satisfying the conditions of corollary 2. Let
${\cal B}_{2}$ and ${\cal B}_{3}$ denote the similar classes
corresponding to corollaries 3 and 4. Let ${\cal B} =
\bigcup_{i=1}^{3}{\cal B}_{i}$.

\vspace{0.05in}\noindent {\sf {\textbf{Theorem 3}}} {\em Let $b_0 :=
1$, $b_k := \sum_{j=1}^{n}w_j b_{j,k}$, $k \geq 1,$ where for every
$j=1, \ldots ,n$, $\{b_{j,k}: k \geq 0\} \in {\cal B}$. Then the
sequence $\{b_k: k \geq 0\}$ satisfies $b_k \geq 1\;\forall\;k \geq
0$, is unbounded and $\{k \geq 0: b_k =1\}$ is finite. Also,
$\lim_{k \rightarrow\infty} b_{k}^{1/k} = 1$. Consequently, $\{b_k:
k \geq 0\}$ satisfies condition (A) of theorem 1 for any $f$ there.}

\vspace{0.05in} \noindent {\sf {\textbf{Proof}}} Notice that we need
only to show that $\lim_{k \rightarrow\infty} b_{k}^{1/k} = 1$. To
see this, note that $1 \leq b_k \leq \max_{1\leq j \leq n} b_{j,k}$.
Therefore, $1 \leq b_k^{1/k} \leq \max_{1\leq j \leq n}
b_{j,k}^{1/k}$. Notice now that by choice of $\{b_{j,k}: k \geq
0\}$, $\lim_{k \rightarrow\infty} b_{j,k}^{1/k} = 1$ for every $j
=1, \ldots ,n$. Hence, $\lim_{k \rightarrow\infty} \max_{1\leq j
\leq n} b_{j,k}^{1/k} = 1$, and this implies, in turn, $\lim_{k
\rightarrow\infty} b_{k}^{1/k} = 1$. \hfill \qed

\vspace{0.05in}\noindent {\sf {\textbf{Remark 6}}} A comparison of
the main result of Bar-Lev and Plachky (1989) for discrete
distributions and our results (theorem 1, theorem 2, corollaries
1--4, and theorem 3) reveals that our construction and results are
different from the one they developed. Also, our results show the
connection with unbiased estimation.

\vspace{0.17in} \noindent {\large\sf {\textbf{4. Some examples for
specific power series distributions}}}

\vspace{0.05in} \noindent We have already seen in section 3
(corollaries 1--4, theorem 3) that there are several choices of
$\left\{b_k: k \geq 0\right\}$ which work for {\bf any} $f$. By way
of illustration of these results, we describe in this section three
new incomplete and boundedly complete families of discrete
distributions. These examples are related to {\bf specific} power
series distributions. We use the notations used in section 1.

In all the examples, we specify $f$ and the sequence $\left\{b_k: k
\geq 0\right\}$. As in theorem 1, ${\cal P}_{f,g} = \{P_{\theta}:
\theta \in (0,R)\}$ denotes the family of discrete distributions
induced by $f$ and $g$. The pmf corresponding to $P_{\theta}$ is
denoted by $p(k;\theta)$. Also, in all the examples, ${\cal
P}_{f,g}$ is an incomplete and boundedly complete family. We specify
$J$, $g(\theta)$, and $\Psi_{0}$ in each.

\vspace{0.05in}\noindent{\sf {\textbf{Example 1}}} ({\em An
illustration of corollary 1 with $n=1$, $w_1 =1$, $d_{1,k} =2$
$\forall\;k \geq 0$}) Let $f(\theta) := e^{\theta},\;\theta \in
{\R}$, and $b_k := 2^k,\;k \geq 0.$ The pmf $p(k;\theta)$
corresponding to $P_{\theta}$ is given by
$$p(k;\theta) = \left\{\begin{array}{lcl} 1-e^{-\theta} & {\mbox{\rm if}} & k = -1, \\
\dpl\frac{e^{-2\theta}\theta^k}{k!} & {\mbox{\rm if}} & k =
0,1,2,\ldots .\end{array}\right.$$

Notice that $J = \{0\}$, $g(\theta)= e^{2\theta}$, $\Psi_{0} =
\{a+c_0 e^{-2\theta}: a, c_0 \in {\R}\}$.

\vspace{0.05in}\noindent{\sf {\textbf{Example 2}}} ({\em An
illustration of corollary 2 with $n=1$, $w_1 =1$}) Let $f(\theta) :=
1-\log(1-\theta),\;\theta \in (0,1)$, $b_0 := 1$, and $b_k := k$ for
$k \geq 1$. The pmf $p(k;\theta)$ corresponding to $P_{\theta}$ is
given by
$$p(k;\theta) = \left\{\begin{array}{lcl} 1-(1-\theta)[1-\log(1-\theta)] & {\mbox{\rm if}} & k = -1, \\
1-\theta & {\mbox{\rm if}} & k = 0,
\\ \dpl\frac{\theta^k(1-\theta)}{k} & {\mbox{\rm if}}
& k = 1,2,\ldots .\end{array}\right.$$

Notice that $J = \{0,1\}$, $g(\theta)= (1-\theta)^{-1}$, $\Psi_{0} =
\{a+c_0(1-\theta)+c_1\theta(1-\theta): a, c_0, c_1 \in {\R}\}$.

\vspace{0.05in} \noindent {\sf {\textbf{Example 3}}} ({\em An
illustration of corollary 3 with $n=1$, $w_1 =1$, $\alpha=1, \beta
=3$}) Let $f(\theta) := (1-\theta)^{-1},\;\theta \in (0,1)$, $b_0 =
1$, and $b_k := t(3,k)/t(1,k) = (k+1)(k+2)/2,\;k \geq 1.$ The pmf
$p(k;\theta)$ corresponding to $P_{\theta}$ is given by
$$p(k;\theta) = \left\{\begin{array}{lcl} 1-(1-\theta)^2 & {\mbox{\rm if}} & k = -1, \\
(1-\theta)^{3}\theta^{k} & {\mbox{\rm if}} & k = 0,1,2,\ldots
.\end{array}\right.$$

Notice that $J = \{0\}$, $g(\theta)= (1-\theta)^{-3}$, and $\Psi_0 =
\{a+c_0 (1-\theta)^{3}: a, c_0 \in {\R}\}$.

\vspace{0.05in}\noindent{\sf {\textbf{Remark 7}}} The example of
Lehmann and Scheffe (1950, p. 312) (or Lehmann and Casella (1998,
pp. 84--85)) is a particular case of corollary 3 with $n=1$, $w_1
=1$, $\alpha=1$, $\beta =2$, where $f(\theta) = (1-\theta)^{-1}$.

\vspace{0.05in}\noindent{\sf {\textbf{Remark 8}}} We believe that
examples 1--3 are simple enough to be taught in an introductory
course in statistical inference or mathematical statistics. The
results in the preceding sections also indicate how one can
construct similar examples.

\vspace{0.17in} \noindent \centerline{\Large\sf {\textbf{Appendix}}}

\vspace{0.05in} \noindent We state below two propositions which were
needed in describing corollary 3. The first one is a standard fact
of elementary analysis.

\vspace{0.05in} \noindent {\sf {\textbf{Proposition 1}}} {\em For
$-1 < \theta < 1$, $\gamma > 0$, $(1-\theta)^{-\gamma} =
\sum_{k=0}^{\infty} t(\gamma,k)\theta^k,$ where}
$$t(\gamma,k) =
\left\{\begin{array}{ll} 1 & {\mbox{\it if}}\;\;k = 0, \\
\dpl\frac{\gamma(\gamma+1)\cdots (\gamma+k-1)}{k!} & {\mbox{\it
if}}\;\;k = 1, 2, \ldots.\end{array}\right.$$

\vspace{0.05in} \noindent {\sf {\textbf{Proposition 2}}} {\em
Suppose $0 < \alpha < \beta$. Then $\lim_{k\rightarrow\infty}
[t(\beta,k)/t(\alpha,k)] = \infty$.}

\noindent {\sf {\textbf{Proof}}} The proof is an immediate
consequence of the following facts: (a) for $k \geq 1,$
$$\dpl\frac{t(\beta,k)}{t(\alpha,k)} = \dpl\prod_{i=1}^{k}\left(1+\dpl\frac{\beta-\alpha}{\alpha+i-1}\right) >
1 + \dpl\sum_{i=1}^{k}\dpl\frac{\beta-\alpha}{\alpha+i-1} > 1 +
\dpl\sum_{i=1}^{k}\dpl\frac{\beta-\alpha}{\lfloor\alpha\rfloor
+i},$$ where $\lfloor \alpha \rfloor$ is the greatest integer less
than or equal to $\alpha$, and (b) $\sum_{n=1}^{\infty}n^{-1} =
\infty$. \hfill \qed

\vspace{0.17in} \noindent \centerline{\Large\sf
{\textbf{References}}}

\vspace{0.15in}\noindent\hang Bar-Lev, S.K. and Plachky, D. (1989).
Boundedly complete families which are not complete. {\em Metrika}
{\bf 36}, 331--336.

\noindent\hang Hoeffding, W. (1977). Some incomplete and boundedly
complete families of distributions. {\em The Annals of Statistics}
{\bf 5}, 278--291.

\noindent\hang Lehmann, E.L. and Casella, G. (1998). {\em Theory of
Point Estimation}, Second Edition. Springer, New York.

\noindent\hang Lehmann, E.L. and Scheffe, H. (1950). Completeness,
similar regions, and unbiased estimation. {\em Sankhy$\bar{a}$} {\bf
10}, 305--340

\noindent\hang Mattner, L. (1993). Some incomplete but boundedly
complete location families. {\em The Annals of Statistics} {\bf 21},
2158--2162.
\end{document}